\documentclass[12pt,twoside]{article}
\usepackage{amsmath, amsthm, amssymb, mathrsfs} 
\usepackage{graphics, graphicx}

\newcommand{\partdif}[2]{\ensuremath{ \frac{\partial #1}{\partial #2}}}

\begin{document}

\title{A Higher Order Godunov Method for Radiation Hydrodynamics: Radiation Subsystem}
\author{Michael Sekora$^{\dagger}$ \\ 
James Stone$^{\dagger,~\ddagger}$ \\ 
\textit{{\small Program in Applied and Computational Mathematics$^{\dagger}$}} \\ 
\textit{{\small Department of Astrophysical Sciences$^{\ddagger}$}} \\ 
\textit{{\small Princeton University, Princeton, NJ 08540, USA}}}
\date{1 December 2008}
\maketitle

\begin{abstract}
\noindent
A higher order Godunov method for the radiation subsystem of radiation hydrodynamics is presented. A key ingredient of the method is the direct coupling of stiff source term effects to the hyperbolic structure of the system of conservation laws; it is composed of a predictor step that is based on Duhamel's principle and a corrector step that is based on Picard iteration. The method is second order accurate in both time and space, unsplit, asymptotically preserving, and uniformly well behaved from the photon free streaming (hyperbolic) limit through the weak equilibrium diffusion (parabolic) limit and to the strong equilibrium diffusion (hyperbolic) limit. Numerical tests demonstrate second order convergence across various parameter regimes. 
\end{abstract}

\pagestyle{myheadings}
\markboth{M. Sekora, J. Stone}{A Higher Order Godunov Method for Radiation Hydrodynamics}


\section{Introduction}
\noindent
Radiation hydrodynamics is a fluid description of matter (plasma) that absorbs and emits electromagnetic radiation and in so doing modifies dynamical behavior. The coupling between matter and radiation is significant in many phenomena related to astrophysics and plasma physics, where radiation comprises a major fraction of the internal energy and momentum and provides the dominant transport mechanism. Radiation hydrodynamics governs the physics of radiation driven outflows, supernovae, accretion disks, and inertial confinement fusion \cite{castorbook, mmbook}. Such physics is described mathematically by a nonlinear system of conservation laws that is obtained by taking moments of the Boltzmann and photon transport equations. A key difficulty is choosing the frame of reference in which to take the moments of the photon transport equation. In the comoving and mixed frame approaches, one captures the matter/radiation coupling by adding relativistic source terms correct to $\mathcal{O}(u/c)$ to the right-hand side of the conservation laws, where $u$ is the material flow speed and $c$ is the speed of light. These source terms are stiff because of the variation in time/length scales associated with such problems \cite{mk1982}. This stiffness causes numerical difficulties and makes conventional methods such as operator splitting and method of lines breakdown \cite{leveque1, leveque2}. \\

\noindent
Previous research in numerically solving radiation hydrodynamical problems was carried out by Castor 1972, Pomraning 1973, Mihalas \& Klein 1982, and Mihalas \& Mihalas 1984 \cite{mmbook, mk1982, castor1972, pomraning1973}. There are a variety of algorithms for radiation hydrodynamics. One of the simplest approaches was developed by Stone, Mihalas, \& Norman 1992 and implemented in the ZEUS code, which was based on operator splitting and Crank-Nicholson finite differencing \cite{stone1992}. Since then, higher order Godunov methods have emerged as a valuable technique for solving hyperbolic conservation laws (e.g., hydrodynamics), particularly when shock capturing and adaptive mesh refinement is important \cite{athena}. However, developing upwind differencing methods for radiation hydrodynamics is a difficult mathematical and computational task. In many cases, Godunov methods for radiation hydrodynamics either: $(i)$ neglect the heterogeneity of weak/strong coupling and solve the system of equations in an extreme limit \cite{dai1, dai2}, $(ii)$ are based on a manufactured limit and solve a new system of equations that attempts to model the full system \cite{jin, buet}, or $(iii)$ uses a variation on flux limited diffusion \cite{levermore1981, heracles2007}. All of these approaches do not treat the full generality of the problem. For example, in a series of papers, Balsara 1999 proposed a Riemann solver for the full system of equations \cite{balsara1999}. However, as pointed out by Lowrie \& Morel 2001, Balsara's method failed to maintain coupling between radiation and matter. Moreover, Lowrie \& Morel were critical of the likelihood of developing a Godunov method for full radiation hydrodynamics \cite{lowrie2001}. \\

\noindent
In radiation hydrodynamics, there are three important dynamical scales and each scale is associated with either the material flow (speed of sound), radiation flow (speed of light), or source terms. When the matter-radiation coupling is strong, the source terms define the fastest scale. However, when the matter-radiation coupling is weak, the source terms define the slowest scale. Given such variation, one aims for a scheme that treats the stiff source terms implicitly. Following work by Miniati \& Colella 2007, this paper presents a method that is a higher order modified Godunov scheme that directly couples stiff source term effects to the hyperbolic structure of the system of conservation laws; it is composed of a predictor step that is based on Duhamel's principle and a corrector step that is based on Picard iteration \cite{mc2007}. The method is explicit on the fastest hyperbolic scale (radiation flow) but is unsplit and fully couples matter and radiation with no approximation made to the full system of equations for radiation hydrodynamics. \\

\noindent
A challenge for the modified Godunov method is its use of explicit time differencing when there is a large range in the time scales associated with the problem, $c / a_{\infty} \gg 1$ where $a_{\infty}$ is the reference material sound speed. One could have built a fully implicit method that advanced time according to the material flow scale, but a fully implicit approach was not pursued because such methods often have difficulties associated with conditioning, are expensive because of matrix manipulation and inversion, and are usually built into central difference schemes rather than higher order Godunov methods. An explicit method may even out perform an implicit method if one considers applications that have flows where $c / a_{\infty} \lesssim 10$. A modified Godunov method that is explicit on the fastest hyperbolic scale (radiation flow) as well as a hybrid method that incorporates a backward Euler upwinding scheme for the radiation components and the modified Godunov scheme for the material components are under construction for full radiation hydrodynamics. A goal of future research is to directly compare these two methods in various limits for different values of $c / a_{\infty}$.


\section{Radiation Hydrodynamics}
\noindent
The full system of equations for radiation hydrodynamics in the Eulerian frame that is correct to $\mathcal{O}(1/\mathbb{C})$ is:
\begin{equation}
\partdif{\rho}{t} + \nabla \cdot \left( \mathbf{m} \right) = 0 , 
\end{equation}
\begin{equation}
\partdif{\mathbf{m}}{t} + \nabla \cdot \left( \frac{\mathbf{m} \otimes \mathbf{m}}{\rho} \right) + \nabla p = -\mathbb{P} \left [ -\sigma_t \left( \mathbf{F_r} - \frac{ \mathbf{u}E_r + \mathbf{u} \cdot \mathsf{P_r} }{\mathbb{C}} \right) + \sigma_a \frac{\mathbf{u}}{\mathbb{C}} (T^4 - E_r) \right ] , 
\end{equation}
\begin{equation}
\partdif{E}{t} + \nabla \cdot \left( (E+p) \frac{\mathbf{m}}{\rho} \right) = -\mathbb{P} \mathbb{C} \left [ \sigma_a(T^4 - E_r) + (\sigma_a - \sigma_s) \frac{\mathbf{u}}{\mathbb{C}} \cdot \left( \mathbf{F_r} - \frac{ \mathbf{u}E_r + \mathbf{u} \cdot \mathsf{P_r} }{\mathbb{C}} \right) \right ] , 
\end{equation}
\begin{equation}
\partdif{E_r}{t} + \mathbb{C} \nabla \cdot \mathbf{F_r} = \mathbb{C} \left [ \sigma_a(T^4 - E_r) + (\sigma_a - \sigma_s) \frac{\mathbf{u}}{\mathbb{C}} \cdot \left( \mathbf{F_r} - \frac{ \mathbf{u}E_r + \mathbf{u} \cdot \mathsf{P_r} }{\mathbb{C}} \right) \right ] \label{eq:rad_sys_e}, 
\end{equation}
\begin{equation}
\partdif{ \mathbf{F_r} }{t} + \mathbb{C} \nabla \cdot \mathsf{P_r} = \mathbb{C} \left [ -\sigma_t \left( \mathbf{F_r} - \frac{ \mathbf{u}E_r + \mathbf{u} \cdot \mathsf{P_r} }{\mathbb{C}} \right) + \sigma_a \frac{\mathbf{u}}{\mathbb{C}} (T^4 - E_r) \right ] \label{eq:rad_sys_f} ,
\end{equation}
\begin{equation}
\mathsf{P_r} = \mathsf{f} E_r ~~ \textrm{(closure relation)} \label{eq:rad_closure} .
\end{equation}

\noindent
For the material quantities, $\rho$ is density, $\mathbf{m}$ is momentum, $p$ is pressure, $E$ is total energy density, and $T$ is temperature. For the radiative quantities, $E_r$ is energy density, $\mathbf{F_r}$ is flux, $\mathsf{P_r}$ is pressure, and $\mathsf{f}$ is the variable tensor Eddington factor. In the source terms, $\sigma_a$ is the absorption cross section, $\sigma_s$ is the scattering cross section, and $\sigma_t = \sigma_a + \sigma_s$ is the total cross section. \\

\noindent
Following the presentation of Lowrie, Morel, \& Hittinger 1999 and Lowrie \& Morel 2001, the above system of equations has been non-dimensionalized with respect to the material flow scale so that one can compare hydrodynamical and radiative effects as well as identify terms that are $\mathcal{O}(u/c)$. This scaling gives two important parameters: $\mathbb{C} = c / a_{\infty}$, $\mathbb{P} = \frac{a_r T^4_{\infty}}{\rho_{\infty}a^2_{\infty}}$. $\mathbb{C}$ measures relativistic effects while $\mathbb{P}$ measures how radiation affects material dynamics and is proportional to the equilibrium radiation pressure over material pressure. $a_r = \frac{8 \pi^5 k^4}{15 c^3h^3}$ is a radiation constant, $T_{\infty}$ is the reference material temperature, and $\rho_{\infty}$ is the reference material density. \\

\noindent
For this system of equations, one has assumed that scattering is isotropic and coherent in the comoving frame, emission is defined by local thermodynamic equilibrium (LTE), and that spectral averages for the cross-sections can be employed (gray approximation). The coupling source terms are given by the modified Mihalas-Klein description \cite{lowrie2001, lowrie1999} which is more general and more accurate than the original Mihalas-Klein source terms \cite{mk1982} because it maintains an important $\mathcal{O}(1 / \mathbb{C}^2)$ term that ensures the correct equilibrium state and relaxation rate to equilibrium \cite{lowrie2001, lowrie1999}. \\

\noindent
Before investigating full radiation hydrodynamics, it is useful to examine the radiation subsystem, which is a simpler system that minimizes complexity while maintaining the rich hyperbolic-parabolic behavior associated with the stiff source term conservation laws. This simpler system allows one to develop a reliable and robust numerical method. Consider Equations \ref{eq:rad_sys_e}, \ref{eq:rad_sys_f} for radiation hydrodynamics in one spatial dimension not affected by transverse flow. If one only considers radiative effects and holds the material flow stationary such that $u \rightarrow 0$, then the conservative variables, fluxes, and source terms for the radiation subsystem are given by:
\begin{equation}
\partdif{E_r}{t} + \mathbb{C}   \partdif{F_r}{x} =  \mathbb{C} \sigma_a (T^4 - E_r) \label{eq:rad_sub_e} ,
\end{equation}
\begin{equation}
\partdif{F_r}{t} + \mathbb{C} f \partdif{E_r}{x} = -\mathbb{C} \sigma_t F_r \label{eq:rad_sub_f} .
\end{equation}

\noindent
Motivated by the asymptotic analysis of Lowrie, Morel, \& Hittinger 1999 for full radiation hydrodynamics, one investigates the limiting behavior for this simpler system of equations. For non-relativistic flows $1/\mathbb{C} = \mathcal{O}(\epsilon)$, where $\epsilon \ll 1$. Assume that there is a moderate amount of radiation in the flow such that $\mathbb{P} = \mathcal{O}(1)$. Furthermore, assume that scattering effects are small such that $\sigma_s / \sigma_t = \mathcal{O}(\epsilon)$. Lastly, assume that the optical depth can be represented as $\mathcal{L} = \ell_{\textrm{mat}} / \lambda_t = \ell_{\textrm{mat}} ~ \sigma_t$, where $\lambda_t$ is the total mean free path of the photos and $\ell_{\textrm{mat}} = \mathcal{O}(1)$ is the material flow length scale \cite{lowrie1999}. \\

\noindent
\textbf{Free Streaming Limit $\sigma_a, \sigma_t \sim \mathcal{O}(\epsilon)$:} In this regime, the right-hand-side of Equations \ref{eq:rad_sub_e} and \ref{eq:rad_sub_f} is negligible such that the system is strictly hyperbolic. $f \rightarrow 1$ and the Jacobian of the quasilinear conservation law has eigenvalues $\pm \mathbb{C}$:
\begin{equation}
\partdif{E_r}{t} + \mathbb{C} \partdif{F_r}{x} = 0 \label{eq:stream_e} , 
\end{equation}
\begin{equation}
\partdif{F_r}{t} + \mathbb{C} \partdif{E_r}{x} = 0 \label{eq:stream_f} .
\end{equation}

\noindent
\textbf{Weak Equilibrium Diffusion Limit $\sigma_a, \sigma_t \sim \mathcal{O}(1)$ :} One obtains this limit by plugging in $\sigma_a, \sigma_t \sim \mathcal{O}(1)$, matching terms of like order, and combining the resulting equations. From the definition of the equilibrium state, $E_r = T^4$ and $F_r = -\frac{1}{\sigma_t} \partdif{P_r}{x}$. Therefore, the system is parabolic and resembles a diffusion equation, where $f \rightarrow 1/3$:
\begin{equation}
\partdif{E_r}{t} = \frac{\mathbb{C}}{3 \sigma_t} \partdif{^2 E_r}{x^2} \label{eq:weak_e} ,
\end{equation}
\begin{equation}
F_r = -\frac{1}{3 \sigma_t} \partdif{E_r}{x} \label{eq:weak_f} .
\end{equation}

\noindent
\textbf{Strong Equilibrium Diffusion Limit $\sigma_a, \sigma_t \sim \mathcal{O}(1/\epsilon)$ :} One obtains this limit by plugging in $\sigma_a, \sigma_t \sim \mathcal{O}(1/\epsilon)$ and following the steps outlined for the weak equilibrium diffusion limit. One can consider the system to be hyperbolic, where $f \rightarrow 1/3$ and the Jacobian of the quasilinear conservation law has eigenvalues $\pm \epsilon$:
\begin{equation}
\partdif{E_r}{t} = 0 \label{eq:strong_e} , 
\end{equation}
\begin{equation}
F_r = 0 \label{eq:strong_f} .
\end{equation}

\noindent
Lowrie, Morel, \& Hittinger 1999 investigated an additional limit for full radiation hydrodynamics, the isothermal regime. This limit has some dynamical properties in common with the weak equilibrium diffusion limit, but its defining characteristic is that the material temperature $T(x,t)$ is constant. When considering the radiation subsystem, there is little difference between the weak equilibrium diffusion and isothermal limits because the material quantities, including the material temperature $T$, do not evolve. $T$ enters the radiation subsystem as a parameter rather than a dynamical quantity.


\section{Higher Order Godunov Method}
\noindent
In one spatial dimension, systems of conservation laws with source terms have the form:
\begin{equation}
\partdif{U}{t} + \partdif{F(U)}{x} = S(U) ,
\end{equation}

\noindent
where $U: \mathbb{R} \times \mathbb{R} \rightarrow \mathbb{R}^n$ is an $n$-dimensional vector of conserved quantities. For the radiation subsystem:
\begin{equation}
U = 
\left( \begin{array}{c} E_r \\
                        F_r \end{array} \right) , ~~
F(U) = 
\left( \begin{array}{c} \mathbb{C} F_r \\
                        \mathbb{C} f E_r \end{array} \right) , ~~
S(U) = 
\left( \begin{array}{c} \mathbb{C} S_E \\
                        \mathbb{C} S_F \end{array} \right) = 
\left( \begin{array}{c}   \mathbb{C} \sigma_a (T^4 - E_r) \\
                        - \mathbb{C} \sigma_t F_r \end{array} \right) . \nonumber
\end{equation}

\noindent 
The quasilinear form of this system of conservation laws is:
\begin{equation}
\partdif{U}{t} + A \partdif{U}{x} = S(U) , ~~ 
A = \partdif{F}{U} = \left( \begin{array}{cc} 
0            & \mathbb{C} \\
\mathbb{C} f & 0 \end{array} \right) .
\end{equation}

\noindent
$A$ has eigenvalues $\lambda = \pm f^{1/2} \mathbb{C}$ as well as right eigenvectors $R$ (stored as columns) and left eigenvectors $L$ (stored as rows):
\begin{equation}
R = \left( \begin{array}{cc} 
1         & 1       \\
- f^{1/2} & f^{1/2} \end{array} \right) , ~~
L = \left( \begin{array}{cc} 
\frac{1}{2} & -\frac{1}{2} \left( \frac{1}{f} \right)^{1/2} \\
\frac{1}{2} &  \frac{1}{2} \left( \frac{1}{f} \right)^{1/2} \end{array} \right) .
\end{equation}

\noindent
Godunov's method obtains solutions to systems of conservation laws by using characteristic information within the framework of a conservative method:
\begin{equation}
U^{n+1}_{i} = U^{n}_{i} - \frac{\Delta t}{\Delta x} \left( F_{i+1/2} - F_{i-1/2} \right) + \Delta t S(U^{n}_{i}) .
\end{equation}

\noindent
Numerical fluxes $F_{i \pm 1/2}$ are obtained by solving the Riemann problem at the cell interfaces with left/right states to get $U_{i-1/2}^{n \pm 1/2}$ and computing $F_{i \pm 1/2} = F(U_{i \pm 1/2}^{n+1/2})$, where $i$ represents the location of a cell center, $i \pm 1/2$ represents the location cell faces to the right and left of $i$, and superscripts represent the time discretization. An HLLE (used in this work) or any other approximate Riemann solver may be employed because the Jacobian $\partial F / \partial U$ for the radiation subsystem is a constant valued matrix and by definition a Roe matrix \cite{leveque1, leveque2, roe1981}. This property also implies that one does not need to transform the system into primitive variables ($\nabla_U W$). The power of the method presented in this paper is that the spatial reconstruction, eigen-analysis, and cell-centered updating directly plug into conventional Godunov machinery. 

\subsection{Predictor Step}
\noindent
One computes the flux divergence $(\nabla \cdot F)^{n+1/2}$ by using the quasilinear form of the system of conservation laws and the evolution along Lagrangian trajectories:
\begin{equation} 
\frac{DU}{Dt} + A^L \partdif{U}{x} = S(U) , ~~ A^L = A - uI, ~~ \frac{DU}{Dt} = \partdif{U}{t} + \left( u \partdif{}{x} \right) U .
\end{equation}

\noindent
From the quasilinear form, one derives a system that includes (at least locally in time and state space) the effects of the stiff source terms on the hyperbolic structure. Following the analysis of Miniati \& Colella 2007 and Trebatich et al 2005 \cite{mc2007, treb2005}, one applies Duhamel's principle to the system of conservation laws, thus giving:
\begin{equation}
\frac{D U^{\textrm{eff}}}{Dt} = \mathcal{I}_{\dot{S}_n}(\eta) \left( - A^L \partdif{U}{x} + S_n \right) ,
\end{equation}

\noindent 
where $\mathcal{I}_{\dot{S}_n}$ is a propagation operator that projects the dynamics of the stiff source terms onto the hyperbolic structure and $\dot{S}_n = \nabla_U S|_{U_n}$. The subscript $n$ designates time $t=t_n$. Since one is considering a first order accurate predictor step in a second order accurate predictor-corrector method, one chooses $\eta = \Delta t / 2$ and the effective conservation law is:
\begin{equation} 
\frac{D U}{Dt} + \mathcal{I}_{\dot{S}_n}(\Delta t / 2) A^L \partdif{U}{x} = \mathcal{I}_{\dot{S}_n}(\Delta t / 2) S_n , ~ \Rightarrow ~ \partdif{U}{t} + A_{\textrm{eff}} \partdif{U}{x} = \mathcal{I}_{\dot{S}_n}(\Delta t / 2) S_n ,
\end{equation}

\noindent 
where $A_{\textrm{eff}} = \mathcal{I}_{\dot{S}_n}(\Delta t / 2) A^L + u I$. In order to compute $\mathcal{I}_{\dot{S}_n}$, one first computes $\dot{S}_n$. Since $\mathbb{C}$, $\sigma_a$, and $\sigma_t$ are constant and one assumes that $\partdif{T}{E_r}, \partdif{T}{F_r} = 0$:
\begin{equation}
\dot{S_n} = \left( \begin{array}{cc} 
-\mathbb{C} \sigma_a & 0 \\
0                    & -\mathbb{C} \sigma_t \end{array} \right) .
\end{equation}

\noindent
$\mathcal{I}_{\dot{S}_n}$ is derived from Duhamel's principle and is given by:
\begin{eqnarray}
\mathcal{I}_{\dot{S_n}} (\Delta t / 2) & = & \frac{1}{\Delta t/2} \int^{\Delta t/2}_{0} e^{\tau \dot{S_n}} d\tau \\
  & = & \left( \begin{array}{cc} 
\alpha & 0     \\
0      & \beta \end{array} \right) , ~~~~ \alpha = \frac{1 - e^{-\mathbb{C} \sigma_a \Delta t / 2}}{\mathbb{C} \sigma_a \Delta t / 2} , ~~ \beta = \frac{1 - e^{-\mathbb{C} \sigma_t \Delta t / 2}}{\mathbb{C} \sigma_t \Delta t / 2} .
\end{eqnarray}

\noindent
Before applying $\mathcal{I}_{\dot{S_n}}$ to $A_L$, it is important to understand that moving-mesh methods can be accommodated in non-relativistic descriptions of radiation hydrodynamics whenever an Eulerian frame treatment is employed. These methods do not require transformation to the comoving frame \cite{lowrie2001}. Since the non-dimensionalization is associated with the hydrodynamic scale, one can use $u_{\textrm{mesh}} = u$ from Lagrangean hydrodynamic methods. \\

\noindent
The effects of the stiff source terms on the hyperbolic structure are accounted for by transforming to a moving-mesh (Lagrangean) frame $A_L = A - uI$, applying the propagation operator $\mathcal{I}_{\dot{S_n}}$ to $A_L$, and transforming back to an Eulerian frame $A_{\textrm{eff}} = \mathcal{I}_{\dot{S_n}} A_L + uI$ \cite{mc2007}. However, because only the radiation subsystem of radiation hydrodynamics is considered $u_{\textrm{mesh}} = u \rightarrow 0$. Therefore, the effective Jacobian is given by:
\begin{equation}
A_{\textrm{eff}} = \left( \begin{array}{cc} 
0                  & \alpha \mathbb{C} \\
\beta f \mathbb{C} & 0                 \end{array} \right) \label{eq:a_eff} ,
\end{equation}

\noindent
which has eigenvalues $\lambda_{\textrm{eff}} = \pm (\alpha \beta)^{1/2} f^{1/2} \mathbb{C}$ with the following limits:
\begin{eqnarray}  
\sigma_a, \sigma_t \rightarrow 0 & \Rightarrow & \alpha, \beta \rightarrow 1 ~ \Rightarrow ~ \lambda_{\textrm{eff}} \rightarrow \pm f^{1/2} \mathbb{C}, ~ (\textrm{free streaming}) \nonumber \\
\sigma_a, \sigma_t \rightarrow \infty & \Rightarrow & \alpha, \beta \rightarrow 0 ~ \Rightarrow ~ \lambda_{\textrm{eff}} \rightarrow \pm \epsilon, ~ (\textrm{strong equilibrium diffusion}) \nonumber .
\end{eqnarray}

\noindent
$A_{\textrm{eff}}$ has right eigenvectors $R_{\textrm{eff}}$ (stored as columns) and left eigenvectors $L_{\textrm{eff}}$ (stored as rows):
\begin{equation}
R_{\textrm{eff}} = \left( \begin{array}{cc} 
1 & 1 \\
-\left ( \frac{\beta f}{\alpha} \right )^{1/2} & \left ( \frac{\beta f}{\alpha} \right )^{1/2} \end{array} \right) , ~~ 
L_{\textrm{eff}} = \left( \begin{array}{cc} 
\frac{1}{2} & -\frac{1}{2} \left ( \frac{\alpha}{\beta f} \right )^{1/2} \\
\frac{1}{2} &  \frac{1}{2} \left ( \frac{\alpha}{\beta f} \right )^{1/2} \end{array} \right) \label{eq:rl_eff} .
\end{equation}

\subsection{Corrector Step}
\noindent
The time discretization for the source term is a single-step, second order accurate scheme based on the ideas from Dutt et al 2000, Minion 2003, and Miniati \& Colella 2007 \cite{mc2007, dutt2000, minion2003}. Given the system of conservation laws, one aims for a scheme that has an explicit approach for the conservative flux divergence term $\nabla \cdot F$ and an implicit approach for the stiff source term $S(U)$. Therefore, one solves a following collection of ordinary differential equations at each grid point:
\begin{equation}
\frac{dU}{dt} = S(U) - ( \nabla \cdot F )^{n+1/2},
\end{equation}

\noindent
where the time-centered flux divergence term is taken to be a constant source which is obtained from the predictor step. Assuming time $t = t_n$, the initial guess for the solution at the next time step is:
\begin{equation}
\hat{U} = U^n + \Delta t (I - \Delta t \nabla_U S(U) |_{U^n})^{-1} (S(U^n) - ( \nabla \cdot F )^{n+1/2}) ,
\end{equation}

\noindent
where:
\begin{equation}
\left( I - \Delta t \nabla_U S(U) \right) = \left( \begin{array}{cc} 
1 + \Delta t \mathbb{C} \sigma_a & 0 \\
0                                & 1 + \Delta t \mathbb{C} \sigma_t \end{array} \right) , 
\end{equation}
\begin{equation}
\left( I - \Delta t \nabla_U S(U) \right)^{-1} = \left( \begin{array}{cc} 
\frac{1}{1 + \Delta t \mathbb{C} \sigma_a} & 0 \\
0                                          & \frac{1}{1 + \Delta t \mathbb{C} \sigma_t} \end{array} \right) .
\end{equation}

\noindent
The error $\epsilon$ is defined as the difference between the initial guess and the solution obtained from the Picard iteration equation where the initial guess was used as a starting value:
\begin{equation}
\epsilon(\Delta t) = U^n + \frac{\Delta t}{2} \left( S(\hat{U}) + S(U^n) \right) - \Delta t ( \nabla \cdot F )^{n+1/2} - \hat{U} .
\end{equation}

\noindent
Following Miniati \& Colella 2007, the correction to the initial guess is given by \cite{mc2007}:
\begin{equation}
\delta(\Delta t) = \left( I - \Delta t \nabla_U S(U) |_{\hat{U}} \right)^{-1} \epsilon(\Delta t) .
\end{equation}

\noindent
Therefore, the solution at time $t = t_n + \Delta t$ is:
\begin{equation}
U^{n+1} = \hat{U} + \delta(\Delta t) \label{eq:update} .
\end{equation}

\subsection{Stability and Algorithmic Issues}
\noindent
The higher order Godunov method satisfies important conditions that are required for numerical stability \cite{mc2007}. First, $\lambda_{\rm{eff}} = \pm (\alpha \beta)^{1/2} f^{1/2} \mathbb{C}$ indicates that the subcharacteristic condition for the characteristic speeds at equilibrium is always satisfied, such that: $\lambda^{-} < \lambda_{\rm{eff}}^{-} < \lambda^0 < \lambda_{\rm{eff}}^{+} < \lambda^{+}$. This condition is necessary for the stability of the system and guarantees that the numerical solution tends to the solution of the equilibrium equation as the relaxation time tends to zero. Second, since the structure of the equations remains consistent with respect to classic Godunov methods, one expects the CFL condition to apply: $\rm{max}(|\lambda^*|) \frac{\Delta t}{\Delta x} \leq 1$, $* = -,0,+$. \\

\noindent
Depending upon how one carries out the spatial reconstruction to solve the Riemann problem in Godunov's method, the solution is either first order accurate in space (piecewise constant reconstruction) or second order accurate in space (piecewise linear reconstruction). Piecewise linear reconstruction was employed in this paper, where left/right states (with respect to the cell center) are modified to account for the stiff source term effects \cite{mc2007, colella1990}:
\begin{eqnarray}
U_{i,\pm}^{n}      & = & U_{i}^{n} + \frac{\Delta t}{2} \mathcal{I}_{\dot{S_n}} \left( \frac{\Delta t}{2} \right) S(U_{i}^{n}) + \frac{1}{2} \left( \pm I - \frac{\Delta t}{\Delta x} A_{\textrm{eff}}^{n} \right) P_{\pm}(\Delta U_i)  \label{eq:spat_recon} \\
P_{\pm}(\Delta U_i) & = & \sum_{\pm \lambda_k > 0} \left( L_{\textrm{eff}}^{k} \cdot \Delta U_i \right) \cdot R_{\textrm{eff}}^{k} \label{eq:slope} .
\end{eqnarray}

\noindent
Left/right one-sided slopes as well as cell center slopes are defined for each cell centered quantity $U_i$. A van Leer limiter is applied to these slopes to ensure monotonicity, thus giving the local slope $\Delta U_i$.


\section{Numerical Tests}
\noindent
Four numerical tests spanning a range of mathematical and physical behavior were carried out to gauge the temporal and spatial accuracy of the higher order Godunov method. The numerical solution is compared with the analytic solution where possible. Otherwise, a self-similar comparison is made. Using piecewise constant reconstruction for the left/right states, one can show that the Godunov method reduces to a consistent discretization in each of the limiting cases. \\

\noindent
The optical depth $\tau$ is a useful quantity for classifying the limiting behavior of a system that is driven by radiation hydrodynamics:
\begin{equation}
\tau = \int_{x_{min}}^{x_{max}} \sigma_t dx = \sigma_t (x_{\max}-x_{\min}) ,
\end{equation}

\noindent
Optically thin/thick regimes are characterized by:
\begin{eqnarray}
\tau &<& O(1) ~~ (\textrm{optically thin}) \nonumber \\
\tau &>& O(1) ~~ (\textrm{optically thick}) \nonumber .
\end{eqnarray}

\noindent
In optically thin regimes (free streaming limit), radiation and hydrodynamics decouple such that the resulting dynamics resembles an advection process. In optically thick regimes (weak/strong equilibrium diffusion limit), radiation and hydrodynamics are strongly coupled and the resulting dynamics resembles a diffusion process. \\

\noindent
The following definitions for the n-norms and convergence rates are used throughout this paper. Given the numerical solution $q^{r}$ at resolution $r$ and the analytic solution $u$, the error at a given point $i$ is: $\epsilon^{r}_{i} = q^{r}_{i} - u$. Likewise, given the numerical solution $q^{r}$ at resolution $r$ and the numerical solution $q^{r+1}$ at the next finer resolution $r+1$ (properly spatially averaged onto the coarser grid), the error resulting from this self-similar comparison at a given point $i$ is: $\epsilon^{r}_{i} = q^{r}_{i} - q^{r+1}_{i}$. The 1-norm and max-norm of the error are:
\begin{equation}
L_1 = \sum_i | \epsilon^{r}_{i} | \Delta x^{r} , ~~~~ L_{\max} = \max_i | \epsilon^{r}_{i} | .
\end{equation}

\noindent
The convergence rate is measured using Richardson extrapolation:
\begin{equation}
R_n = \frac{ \textrm{ln} \left( L_n(\epsilon^r)/L_n(\epsilon^{r+1}) \right) }{ \textrm{ln} \left( \Delta x^{r}/\Delta x^{r+1} \right) } .
\end{equation}

\subsection{Exponential Growth/Decay to Thermal Equilibrium}
\noindent
The first numerical test examines the temporal accuracy of how variables are updated in the corrector step. Given the radiation subsystem and the following initial conditions:
\begin{displaymath}
E_r^0 = \textrm{constant across space} , ~~
F_r^0 = 0, ~~
T     = \textrm{constant across space} ,
\end{displaymath}

\noindent
$F_r \rightarrow 0$ for all time. Therefore, the radiation subsystem reduces to the following ordinary differential equation:
\begin{equation}
\frac{d E_r}{dt} = \mathbb{C} \sigma_a (T^4-E_r),
\end{equation}

\noindent 
which has the following analytic solution:
\begin{equation}
E_r = T^4 + (E_r^0 - T^4) \rm{exp}(-\mathbb{C} \sigma_a t) .
\end{equation}

\noindent
For $E_r^0 < T^4$ and $F_r^0 = 0$, one expects exponential growth in $E_r$ until thermal equilibrium $(E_R = T^4)$ is reached. For $E_r^0 > T^4$ and $F_r^0 = 0$, one expects exponential decay in $E_r$ until thermal equilibrium is reached. This numerical test allows one to examine the order of accuracy of the stiff ODE integrator. \\

\noindent
\textbf{Parameters:}
\begin{equation}
\mathbb{C} = 10^{5} , ~ \sigma_{a} = 1 , ~ \sigma_{t} = 2 , ~ f = 1 \nonumber ,
\end{equation}
\begin{equation}
N_{cell} = [32, ~ 64, ~ 128, ~ 256] \nonumber , 
\end{equation}
\begin{equation}
x_{\min} = 0 , ~ x_{\max} = 1 , ~ \Delta x  = \frac{x_{\min}-x_{\max}}{N_{cell}} , ~ CFL = 0.5 , ~ \Delta t  = \frac{CFL ~ \Delta x}{f^{1/2} \mathbb{C}} \nonumber ,
\end{equation}
\begin{equation}
\textrm{IC for Growth:} ~~ E_r^0 = 1, ~ F_r^0 = 0 , ~ T = 10 \nonumber ,
\end{equation}
\begin{equation}
\textrm{IC for Decay:} ~~ E_r^0 = 10^4, ~ F_r^0 = 0 , ~ T = 1 \nonumber .
\end{equation}

\begin{figure}[h]
\begin{center}
\begin{minipage}{3in}
\includegraphics[width=3in,angle=0]{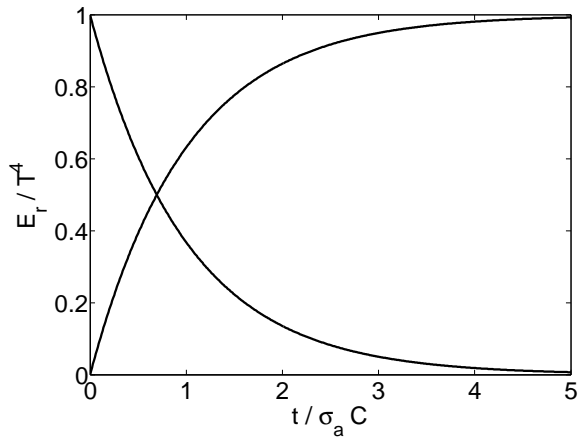}
\caption{\label{fig:exp}Exponential growth/decay to thermal equilibrium. $N_{cell}=256$.}
\end{minipage}
\end{center}
\end{figure} 

\begin{table}[h]
\begin{center}
\begin{tabular*}{1.0\textwidth}{@{\extracolsep{\fill}}rcccccccc}
\hline
$N_{cell}$ & $L_1(E^{g}_{r})$ & Rate & $L_{\infty}(E^{g}_{r})$ & Rate & $L_1(E^{d}_{r})$ & Rate & $L_{\infty}(E^{g}_{r})$ & Rate \\
\hline
 32 & 1.4E-1 & -   & 1.4E-1 & -   & 1.4E-1 & -   & 1.4E-1 & -   \\
 64 & 3.7E-2 & 2.0 & 3.7E-2 & 2.0 & 3.7E-2 & 2.0 & 3.7E-2 & 2.0 \\
128 & 9.3E-3 & 2.0 & 9.3E-3 & 2.0 & 9.3E-3 & 2.0 & 9.3E-3 & 2.0 \\
256 & 2.3E-3 & 2.0 & 2.3E-3 & 2.0 & 2.3E-3 & 2.0 & 2.3E-3 & 2.0 \\
\hline
\label{tbl:exp}
\end{tabular*}
Table 1: Errors and convergence rates for exponential growth/decay in $E_r$ to thermal equilibrium. Errors were obtained through analytic comparison. $t = 10^{-5} = 1/\sigma_a \mathbb{C}$.
\end{center} 
\end{table}

\noindent
From Figure \ref{fig:exp}, one sees that the numerical solution corresponds with the analytic solution. In Table 1, the errors and convergence rates are identical for growth and decay. This symmetry illustrates the robustness of the Godunov method. Furthermore, one finds that the method is well behaved and obtains the correct solution with second order accuracy for stiff values of the $e$ folding time $(\frac{\Delta t}{1 / \sigma_{a} \mathbb{C}} \geq 1)$, although with a significantly larger amplitude in the norm of the error. This result credits the flexibility of the temporal integrator in the corrector step. \\

\noindent
In a similar test, the initial conditions for the radiation energy and flux are zero and the temperature is defined by some spatially varying profile (a Gaussian pulse). As time increases, the radiation energy grows into $T(x)^4$. Unless the opacity is sufficiently high, the radiation energy approaches but does not equal $T(x)^4$. This result shows that the solution has reached thermal equilibrium and any spatially varying temperature will diffuse.

\subsection{Free Streaming Limit}
\noindent
In the free streaming limit, $\tau \ll O(1)$ and the radiation subsystem reduces to Equations \ref{eq:stream_e}, \ref{eq:stream_f}. If one takes an additional temporal and spatial partial derivative of the radiation subsystem in the free streaming limit and subtracts the resulting equations, then one finds two decoupled wave equations that have the following analytic solutions:
\begin{eqnarray}
E_r(x,t) = E_0(x - f^{1/2} \mathbb{C} t) , \\
F_r(x,t) = F_0(x - f^{1/2} \mathbb{C} t) .
\end{eqnarray}

\noindent
\textbf{Parameters:}
\begin{equation}
\mathbb{C} = 10^{5} , ~ \sigma_{a} = 10^{-6} , ~ \sigma_{t} = 10^{-6} , ~ f = 1 , ~ T = 1 , \nonumber 
\end{equation}
\begin{equation}
N_{cell} = [32, ~ 64, ~ 128, ~ 256] , \nonumber 
\end{equation}
\begin{equation}
x_{\min} = 0 , ~ x_{\max} = 1 , ~ \Delta x  = \frac{x_{\min}-x_{\max}}{N_{cell}} , ~ CFL = 0.5 , ~ \Delta t  = \frac{CFL ~ \Delta x}{f^{1/2} \mathbb{C}} , \nonumber
\end{equation}
\begin{equation}
\textrm{IC for Gaussian Pulse:} ~~ E_r^0, F_r^0 = \exp \left( -(\nu (x - \mu) )^2 \right) , ~ \nu = 20 , ~ \mu = 0.3 , \nonumber
\end{equation}
\begin{equation}
\textrm{IC for Square Pulse:} ~~ 
E_r^0, F_r^0 = \left\{ \begin{array}{ll}
               1 & 0.2 < x < 0.4 \\
               0 & \rm{otherwise} \end{array} \right . \nonumber
\end{equation}

\begin{figure}
\begin{center}
\begin{minipage}{3in}
\includegraphics[width=3in,angle=0]{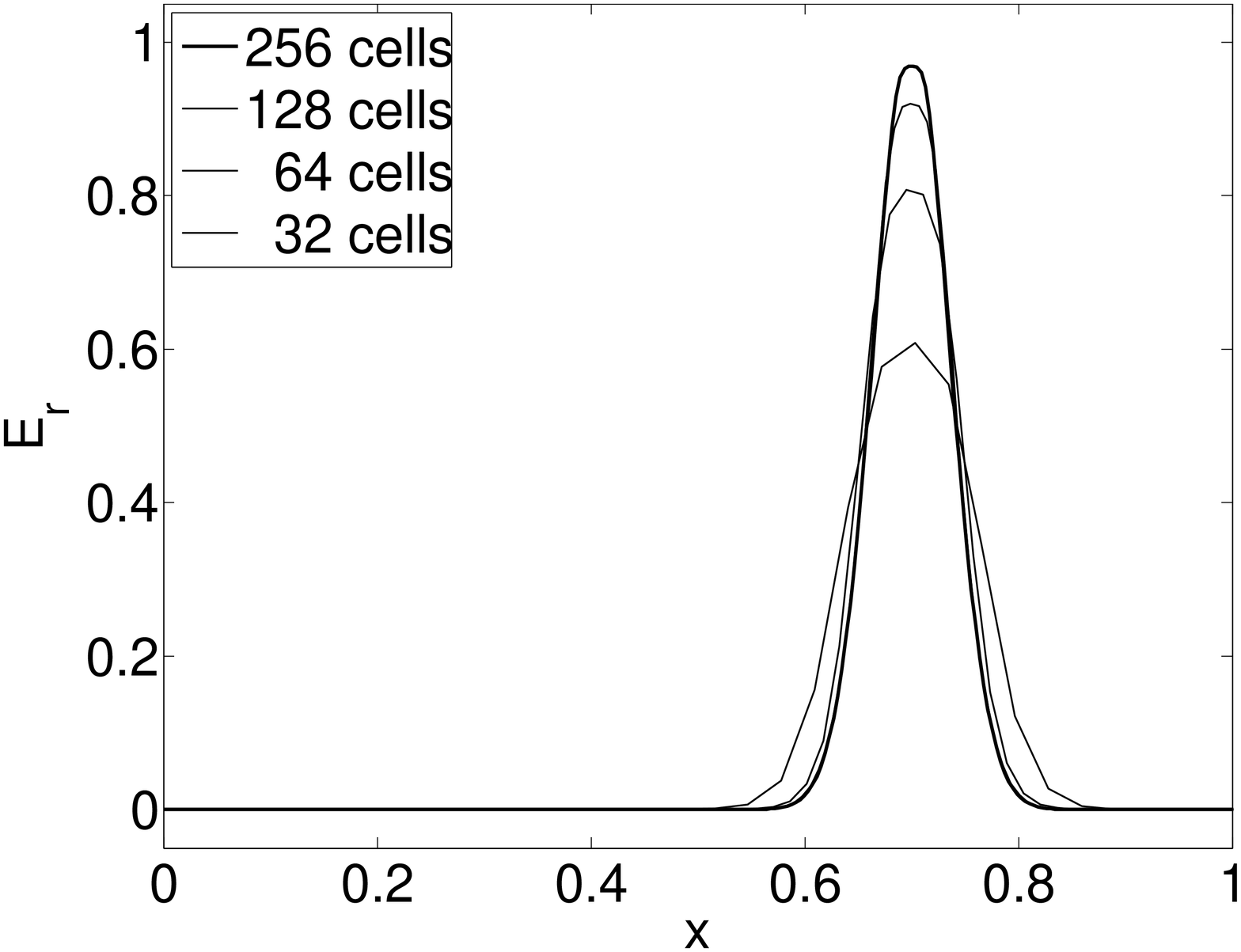}
\caption{\label{fig:stream_gauss} Gaussian pulse in free streaming limit. $t = 4 \times 10^{-6} = 0.4 ~ (x_{\max}-x_{\min}) / \mathbb{C}$.}
\end{minipage} 
\hspace{0.4in}
\begin{minipage}{3in}
\includegraphics[width=3in,angle=0]{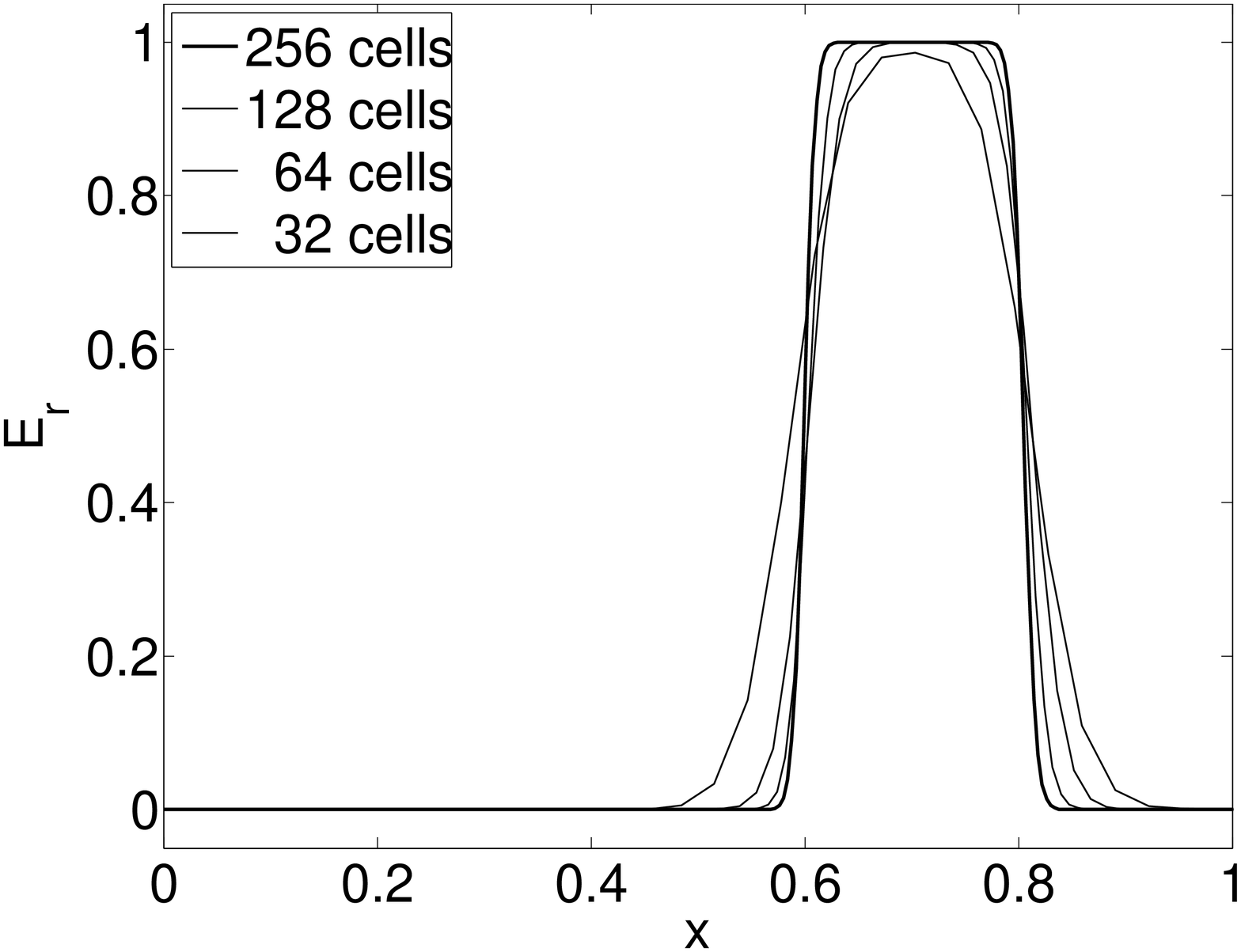}
\caption{\label{fig:stream_sq} Square pulse in free streaming limit. $t = 4 \times 10^{-6} = 0.4 ~ (x_{\max}-x_{\min}) / \mathbb{C}$.}
\end{minipage}
\end{center}
\end{figure} 

\begin{table}
\begin{center}
\begin{tabular*}{1.0\textwidth}{@{\extracolsep{\fill}}rcccccccc}
\hline
$N_{cell}$ & $L_1(E_{r})$ & Rate & $L_{\infty}(E_{r})$ & Rate & $L_1(F_{r})$ & Rate & $L_{\infty}(F_{r})$ & Rate \\
\hline
 32 & 3.8E-2 & -   & 3.9E-1 & -   & 3.8E-2 & -   & 3.9E-1 & -   \\
 64 & 1.3E-2 & 1.5 & 1.8E-1 & 1.1 & 1.3E-2 & 1.5 & 1.8E-1 & 1.1 \\
128 & 3.6E-3 & 1.9 & 8.0E-2 & 1.2 & 3.6E-3 & 1.9 & 8.0E-2 & 1.2 \\
256 & 8.6E-4 & 2.1 & 3.1E-2 & 1.4 & 8.6E-4 & 2.1 & 3.1E-2 & 1.4 \\
\hline
\label{tbl:stream_gauss}
\end{tabular*}
Table 2: Errors and convergence rates for Gaussian pulse in free streaming limit. Errors were obtained through analytic comparison. $t = 4 \times 10^{-6} = 0.4 ~ (x_{\max}-x_{\min}) / \mathbb{C}$.
\end{center} 
\end{table}

\begin{table}
\begin{center}
\begin{tabular*}{1.0\textwidth}{@{\extracolsep{\fill}}rcccc}
\hline
$N_{cell}$ & $L_1(E_{r})$ & Rate & $L_1(F_{r})$ & Rate \\
\hline
 32 & 6.0E-2 & -   & 6.0E-2 & -   \\
 64 & 4.2E-2 & 0.5 & 4.2E-2 & 0.5 \\
128 & 2.6E-2 & 0.7 & 2.6E-2 & 0.7 \\
256 & 1.5E-2 & 0.8 & 1.5E-2 & 0.8 \\
\hline
\label{tbl:stream_sq}
\end{tabular*}
Table 3: Errors and convergence rates for square pulse in free streaming limit. Errors were obtained through analytic comparison. $t = 4 \times 10^{-6} = 0.4 ~ (x_{\max}-x_{\min}) / \mathbb{C}$.
\end{center} 
\end{table}

\noindent
Since the Gaussian pulse results from smooth initial data, one expects $R_1=2.0$. However, the square wave results from discontinuous initial data and one expects $R_1 \simeq 0.67$. This claim is true for all second order spatially accurate numerical methods when applied to an advection-type problem $(u_t + a u_x = 0)$ \cite{leveque1}.

\subsection{Weak Equilibrium Diffusion Limit} 
\noindent
In the weak equilibrium diffusion limit, $\tau > O(1)$ and the radiation subsystem reduces to Equations \ref{eq:weak_e}, \ref{eq:weak_f}. The optical depth suggests the range of total opacities for which diffusion is observed: if $\tau = \sigma_t ~ \ell_{\textrm{diff}} > 1$, then one expects diffusive behavior for $\sigma_t > 1 / \ell_{\textrm{diff}}$. Additionally, Equations \ref{eq:weak_e}, \ref{eq:weak_f} set the time scale $t_{\textrm{diff}}$ and length scale $\ell_{\textrm{diff}}$ for diffusion, where $t_{\textrm{diff}} \sim \ell_{\textrm{diff}}^{~2}/D$ and $D = f \mathbb{C} / \sigma_{t}$ for the radiation subsystem. Given a diffusion problem for a Gaussian pulse defined over the entire real line $(u_t - D u_{xx} = 0)$, the analytic solution is given by the method of Green's functions:
\begin{equation}
u(x,t) = \int_{-\infty}^{\infty} f(\bar{x}) G(x,t;\bar{x},0) d\bar{x} \nonumber = \frac{1}{(4 D t \nu^2 + 1)^{1/2}} \textrm{exp} \left( \frac{-(\nu (x - \mu) )^2}{4 D t \nu^2 + 1} \right) . \label{diff_anal}
\end{equation}

\noindent
\textbf{Parameters:}
\begin{equation}
\mathbb{C} = 10^{5} , ~ \sigma_{a} = 40 , ~ \sigma_{t} = 40 , ~ f = 1/3 , ~ T^4 = E_r , \nonumber
\end{equation}
\begin{equation}
N_{cell} = [320, ~ 640, ~ 1280, ~ 2560] , \nonumber
\end{equation}
\begin{equation}
x_{\min} = -5 , ~ x_{\max} = 5 , ~ \Delta x  = \frac{x_{\min}-x_{\max}}{N_{cell}} , ~ CFL = 0.5 , ~ \Delta t  = \frac{CFL ~ \Delta x}{f^{1/2} \mathbb{C}} , \nonumber
\end{equation}
\begin{equation}
\textrm{IC for Gaussian Pulse:} ~~ \left\{ \begin{array}{ll}
E_r^0 = \exp \left( -(\nu (x - \mu) )^2 \right) , ~ \nu = 20 , ~ \mu = 0.3 , \\
F_r^0 = -\frac{f}{\sigma_t} \partdif{E_r^0}{x} = \frac{2 f \nu^2 (x-\mu)}{\sigma_t} E_r^0 \end{array} \right . \nonumber
\end{equation}

\begin{figure}
\begin{center}
\begin{minipage}{3in}
\includegraphics[width=3in,angle=0]{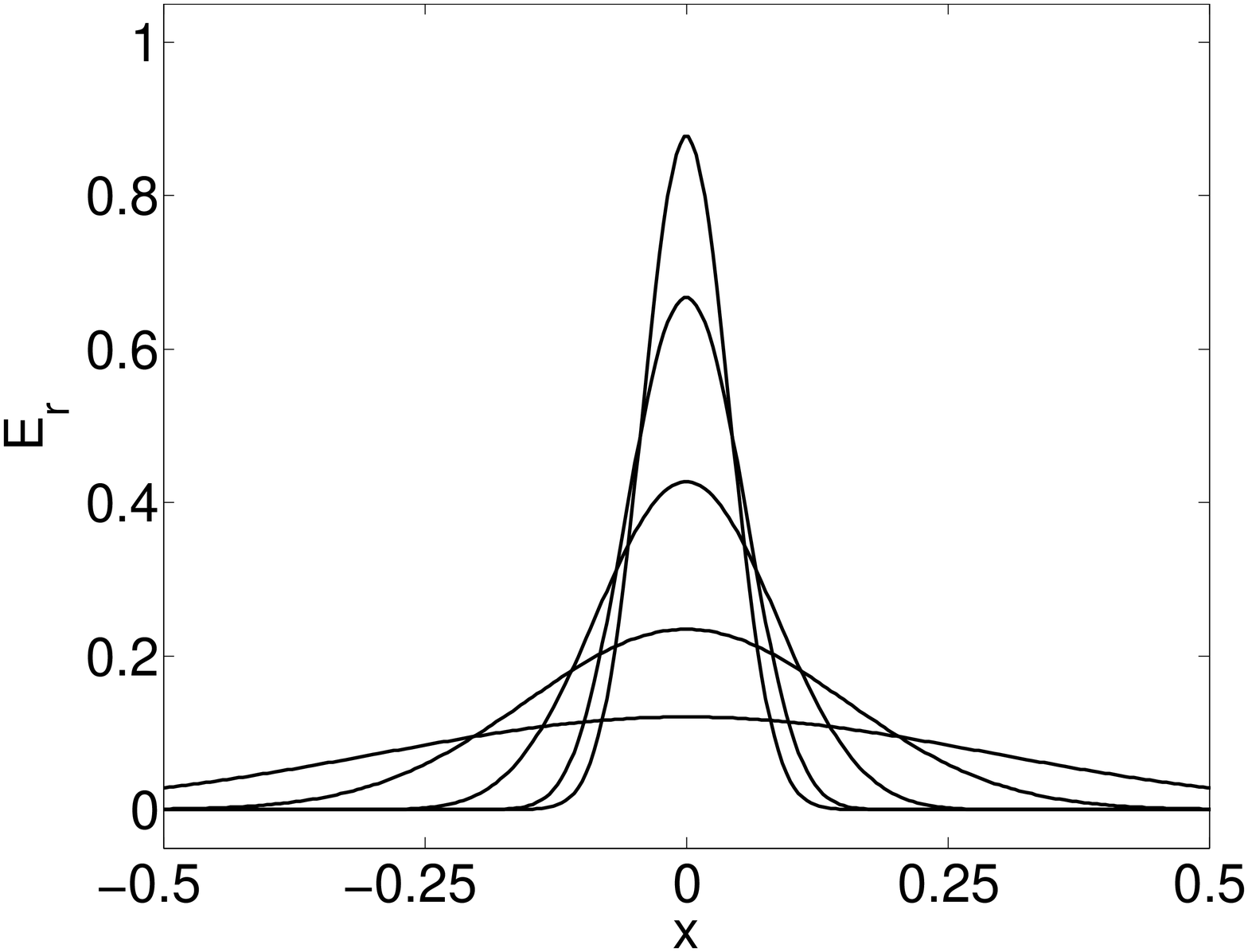}
\caption{\label{fig:weak_e} $E_r$ in weak equilibrium diffusion limit. $t = [0.25, ~ 1, ~ 4, ~ 16, ~ 64] \times 10^{-6}$.}
\end{minipage} 
\hspace{0.4in}
\begin{minipage}{3in}
\includegraphics[width=3in,angle=0]{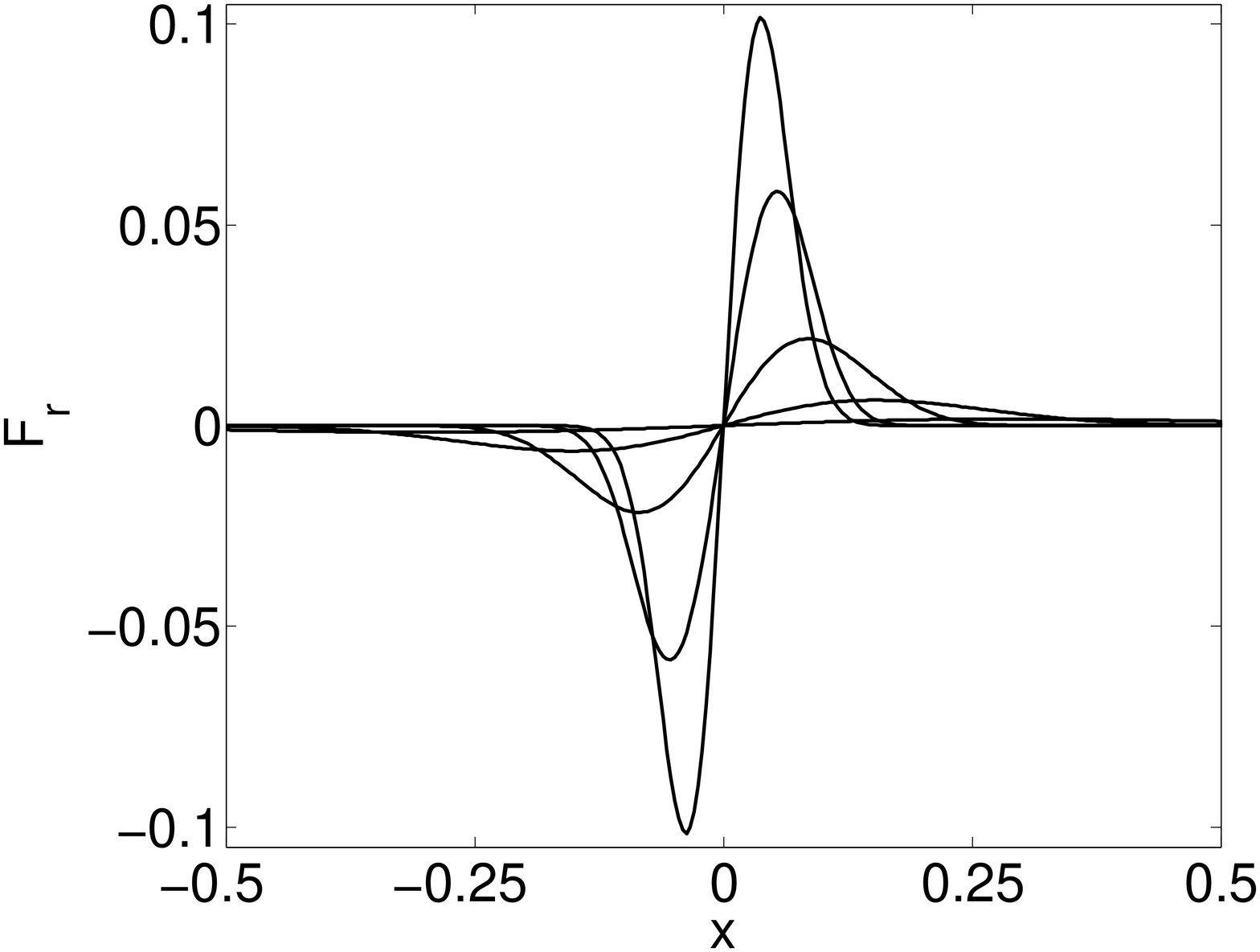}
\caption{\label{fig:weak_f} $F_r$ in weak equilibrium diffusion limit. $t = [0.25, ~ 1, ~ 4, ~ 16, ~ 64] \times 10^{-6}$.}
\end{minipage}
\end{center}
\end{figure} 

\begin{table}
\begin{center}
\begin{tabular*}{1.0\textwidth}{@{\extracolsep{\fill}}rcccccccc}
\hline
$N_{cell}$ & $L_1(E_{r})$ & Rate & $L_{\infty}(E_{r})$ & Rate & $L_1(F_{r})$ & Rate & $L_{\infty}(F_{r})$ & Rate \\
\hline
 320 & 8.9E-3 & -   & 4.5E-2 & -   & 1.1E-3 & -   & 3.7E-3 & -   \\
 640 & 6.6E-3 & 0.4 & 3.4E-2 & 0.4 & 8.3E-4 & 0.4 & 3.1E-3 & 0.2 \\
1280 & 3.4E-3 & 1.0 & 1.6E-2 & 1.1 & 4.1E-4 & 1.0 & 1.4E-3 & 1.2 \\
2560 & 1.6E-3 & 1.1 & 7.1E-3 & 1.1 & 1.9E-4 & 1.1 & 6.0E-4 & 1.2 \\
\hline
\label{tbl:weak_th}
\end{tabular*}
Table 4: Errors and convergence rates for $E_r$, $F_r$ in the weak equilibrium diffusion limit. Time was advanced according to a hyperbolic time step: $\Delta t_{h} = \frac{CFL ~ \Delta x}{f^{1/2} \mathbb{C}}$. Errors were obtained through self-similar comparison. $t = 4 \times 10^{-6}$.
\end{center}
\end{table}

\begin{table}
\begin{center}
\begin{tabular*}{1.0\textwidth}{@{\extracolsep{\fill}}rcccccccc}
\hline
$N_{cell}$ & $L_1(E_{r})$ & Rate & $L_{\infty}(E_{r})$ & Rate & $L_1(F_{r})$ & Rate & $L_{\infty}(F_{r})$ & Rate \\
\hline
 320 & 1.7E-2 & -   & 8.3E-2 & -   & 2.0E-3 & -   & 7.9E-3 & -   \\
 640 & 5.0E-3 & 1.7 & 2.5E-2 & 1.7 & 6.0E-4 & 1.7 & 2.0E-3 & 2.0 \\
1280 & 1.1E-3 & 2.2 & 5.1E-3 & 2.3 & 1.3E-4 & 2.3 & 3.6E-4 & 2.4 \\
2560 & 2.5E-4 & 2.1 & 1.2E-3 & 2.1 & 2.8E-5 & 2.2 & 7.4E-5 & 2.3 \\
\hline
\label{tbl:weak_tp}
\end{tabular*}
Table 5: Errors and convergence rates for $E_r$, $F_r$ in the weak equilibrium diffusion limit. Time was advanced according to a parabolic time step: $\Delta t_{p} = \frac{CFL ~ (\Delta x)^2}{2 D}$. Errors were obtained through self-similar comparison. $t = 4 \times 10^{-6}$.
\end{center} 
\end{table}

\noindent
One's intuition about diffusive processes is based on considering an infinite domain. So to minimize boundary effects in the numerical calculation, the computational domain and number of grid cells were expanded by a factor of 10. In Figures \ref{fig:weak_e}, \ref{fig:weak_f}, one observes the diffusive behavior expected for this parameter regime. Additionally, the numerical solution compares well with the analytic solution for a diffusion process defined over the entire real line (Equation \ref{diff_anal}). However, diffusive behavior is only a first order approximation to more complicated hyperbolic-parabolic dynamics taking place in radiation hydrodynamics as well as the radiation subsystem. Therefore, one needs to compare the numerical solution self-similarly. In Table 4, one sees first order convergence when a hyperbolic time step $\Delta t_{h} = \frac{CFL ~ \Delta x}{f^{1/2} \mathbb{C}}$ is used; while in Table 5, one sees second order convergence when a parabolic time step $\Delta t_{p} = \frac{CFL ~ (\Delta x)^2}{2 D}$ is used. This difference in the convergence rate results from the temporal accuracy in the numerical solution. In the weak equilibrium diffusion limit, the Godunov method reduces to a forward-time/centered-space discretization of the diffusion equation. Such a discretization requires a parabolic time step $\Delta t \sim (\Delta x)^2$ in order to see second order convergence because the truncation error of the forward-time/centered-space discretization of the diffusion equation is $\mathcal{O}(\Delta t,(\Delta x)^2)$. 

\subsection{Strong Equilibrium Diffusion Limit}
\noindent
In the strong equilibrium diffusion limit, $\tau \gg O(1)$. From Equations \ref{eq:strong_e}, \ref{eq:strong_f}, $F_r \rightarrow 0$ for all time and space while $E_r = E_r^0$. \\

\noindent
\textbf{Parameters:}
\begin{equation}
\mathbb{C} = 10^{5} , ~ \sigma_{a} = 10^6 , ~ \sigma_{t} = 10^6 , ~ f = 1/3 , ~ T^4 = E_r , \nonumber
\end{equation}
\begin{equation}
N_{cell} = [320, ~ 640, ~ 1280, ~ 2560] , \nonumber
\end{equation}
\begin{equation}
x_{\min} = -5 , ~ x_{\max} = 5 , ~ \Delta x  = \frac{x_{\min}-x_{\max}}{N_{cell}} , ~ CFL = 0.5 , ~ \Delta t  = \frac{CFL ~ \Delta x}{f^{1/2} \mathbb{C}} , \nonumber
\end{equation}
\begin{equation}
\textrm{IC for Gaussian Pulse:} ~~ \left\{ \begin{array}{ll}
E_r^0 = \exp \left( -(\nu (x - \mu) )^2 \right) , ~ \nu = 20 , ~ \mu = 0.3 , \\
F_r^0 = -\frac{f}{\sigma_t} \partdif{E_r^0}{x} = \frac{2 f \nu^2 (x-\mu)}{\sigma_t} E_r^0 \end{array} \right . \nonumber
\end{equation}

\begin{table}
\begin{center}
\begin{tabular*}{1.0\textwidth}{@{\extracolsep{\fill}}rcccc}
\hline
$N_{cell}$ & $L_1(E_{r})$ & Rate & $L_{\infty}(E_{r})$ & Rate \\
\hline
 320 & 2.2E-3 & -   & 1.8E-2 & -   \\
 640 & 5.3E-4 & 2.1 & 5.6E-3 & 1.6 \\
1280 & 1.3E-4 & 2.0 & 1.5E-3 & 1.9 \\
2560 & 3.3E-5 & 2.0 & 3.8E-4 & 2.0 \\
\hline
\label{tbl:strong}
\end{tabular*}
Table 6: Errors and convergence rates for $E_r$ in the strong equilibrium diffusion limit. Errors were obtained through self-similar comparison. $t = 4 \times 10^{-6}$.
\end{center} 
\end{table}

\noindent
In this test, the numerical solution is held fixed at the initial distribution because $\sigma_a, \sigma_t$ are so large. However, if one fixed $\ell_{\textrm{diff}}$ and scaled time according to $t_{\textrm{diff}} \approx \ell_{\textrm{diff}}^{~2}/D = \ell_{\textrm{diff}}^{~2} \sigma_t / f \mathbb{C}$, then one would observe behavior similar to Figures \ref{fig:weak_e}, \ref{fig:weak_f}. This test illustrates the robustness of the Godunov method to handle very stiff source terms.


\section{Conclusions and Future Work}
\noindent
This paper presents a Godunov method for the radiation subsystem of radiation hydrodynamics that is second order accurate in both time and space, unsplit, asymptotically preserving, and uniformly well behaved. Moreover, the method employs familiar algorithmic machinery without a significant increase in computational cost. This work is the starting point for developing a Godunov method for full radiation hydrodynamics. The ideas in this paper should easily extend to the full system in one and multiple dimensions using a MUSCL or CTU approach \cite{colella1990}. A modified Godunov method that is explicit on the fastest hyperbolic scale (radiation flow) as well as a hybrid method that incorporates a backward Euler upwinding scheme for the radiation components and the modified Godunov scheme for the material components are under construction for full radiation hydrodynamics. A goal of future research is to directly compare these two methods in various limits for different values of $c / a_{\infty}$. Nevertheless, one expects the modified Godunov method that is explicit on the fastest hyperbolic scale to exhibit second order accuracy for all conservative variables and the hybrid method to exhibit first order accuracy in the radiation variables and second order accuracy in the material variables. Work is also being conducted on applying short characteristic and Monte Carlo methods to solve the photon transport equation and obtain the variable tensor Eddington factors. In the present work, these factors were taken to be constant in their respective limits.


\section*{Acknowledgment}
\noindent
The authors thank Dr. Phillip Colella for many helpful discussions. MS acknowledges support from the DOE CSGF Program which is provided under grant DE-FG02-97ER25308. JS acknowledges support from grant DE-FG52-06NA26217.



\end{document}